\documentclass[12pt,oneside]{amsart} 
\usepackage{amssymb,amscd}
\usepackage{euscript}

      \theoremstyle{plain}
      \newtheorem{theorem}{Theorem}[section]
      \newtheorem{lemma}[theorem]{Lemma}
      \newtheorem{corollary}[theorem]{Corollary}
      \newtheorem{proposition}[theorem]{Proposition}
      \newtheorem{remark}[theorem]{Remark}

\numberwithin{equation}{section}

      \makeatletter
      \def\@setcopyright{}
      \def\serieslogo@{}
      \makeatother

\def\B{\EuScript{B}} 
\def\M{\EuScript{M}}
\def\c{\EuScript{C}}
\def\R{\mathbb R}
\def\C{\mathbb C}
\def\Z{\mathbb Z}
\def\N{\mathbb N}
\def\T{\mathbb T}

\def\dist{\text{dist}}

\def\Id{\text{Id}}
\def\e{\varepsilon}
\def\a{\alpha}

\def\tt{\tau}

\def\bWu{\bar{W}^u}
\def\bWs{\bar{W}^s}

\def\bEu{\bar{E}^u}
\def\bEs{\bar{E}^s}
\def\bEc{\bar{E}^c}
\def\f{\bar f}
\def\bm{\bar \mu}
\def\bM{\bar {\EuScript{M}}}

\def\QED{\hfill\hfill{\square}}

\begin{document}

\date{\today}
\author{Boris Kalinin$^\ast$ and Victoria Sadovskaya$^{\ast\ast}$}

\address{Department of Mathematics $\&$ Statistics, 
 University of South  Alabama, Mobile, AL 36688, USA}
\email{kalinin@jaguar1.usouthal.edu, sadovska@jaguar1.usouthal.edu}

\title [Anosov diffeomorphisms  with asymptotically conformal periodic data $\;\;$]
{On Anosov diffeomorphisms  with asymptotically conformal periodic data} 

\thanks{$^{\ast}$  Supported in part by NSF grant DMS-0701292}
\thanks{$^{\ast\ast}$ Supported in part by NSF grant DMS-0401014}


\begin{abstract}
We consider transitive Anosov diffeomorphisms for which every periodic
orbit has only one positive and one negative Lyapunov exponent.  We
establish various properties of such systems including strong pinching, 
$C^{1+\beta}$ smoothness of the Anosov splitting, and $C^1$ smoothness
of measurable invariant conformal structures and distributions.
We apply these results to volume preserving diffeomorphisms with 
two-dimensional stable and unstable distributions and diagonalizable derivatives 
of the return maps at periodic points.  We show that a finite cover of such 
a diffeomorphism is smoothly conjugate to an Anosov automorphism of $\T^4$.
As a corollary we obtain local rigidity for such diffeomorphisms. 
We also establish a local rigidity result for Anosov diffeomorphisms in dimension three.    
\end{abstract}

\maketitle 


\section{Introduction}

The goal of this paper is to study Anosov diffeomorphisms for which every 
periodic orbit has only one positive and one negative Lyapunov exponent.
Our main motivation comes from the problem of local rigidity for higher-dimensional
Anosov systems, i.e. the question of regularity of conjugacy to a small 
perturbation. If $f$ is an Anosov diffeomorphism and $g$ is sufficiently 
$C^1$ close to $f$, then it is well known that $g$ is also Anosov and  
topologically conjugate to $f$. However, the conjugacy is typically 
only H\"older continuous. A necessary condition
for the conjugacy to be $C^1$ is that Jordan normal forms of the derivatives 
of the return maps of $f$ and $g$ at the corresponding periodic points are 
the same. If this condition is also sufficient for any $g$ which is $C^1$ close 
to $f$, then $f$ is called {\em locally rigid}. The problem of local rigidity has 
been extensively studied, and Anosov diffeomorphisms with one-dimensional 
stable and unstable distributions were shown to be locally rigid \cite{L0},
\cite{L1}, \cite{LM}, \cite{P}. 

In contrast, higher dimensional systems are not always locally rigid. In \cite{L1, L2} 
R. de la Llave constructed examples of Anosov automorphisms of the torus $\T^4$
which are not $C^1$ conjugate to certain small perturbations with the same
periodic data. One of the examples has two (un)stable eigenvalues of 
different moduli and the other one has a double (un)stable eigenvalue with a 
nontrivial Jordan block. This suggests that it is natural to consider  local rigidity 
 for  automorphisms that are diagonalizable over $\C$ with all 
(un)stable eigenvalues equal in modulus. For such 
an automorphism the expansion of the unstable and contraction of the stable 
distribution are {\em conformal} with respect to some metric on the torus. For the case of an Anosov diffeomorphism of a compact manifold,
a more natural notion is that of {\em uniform quasiconformality}.
It means that at each point all vectors in the (un)stable subspace 
are expanded/contracted at essentially the same rate (see Section \ref{prelim qc}). 

The study of local rigidity of conformal and uniformly quasiconformal
Anosov systems was initiated in \cite{L2} and continued in \cite{KaS}
and \cite{L4}. It is closely related to the study of {\em global rigidity} of such 
systems, or their classification up to a smooth conjugacy. In \cite{KaS} we 
established the following global rigidity result:

\vskip.2cm
\cite[Theorem 1.1]{KaS} {\it    Let $f$ be a transitive $C^\infty$ Anosov 
   diffeomorphism of a 
   compact manifold $\M$ which is uniformly quasiconformal on the  
   stable and unstable distributions.
   Suppose either that both distributions have dimension at least three,
   or that they have dimension at least two and $\M$ is an
   infranilmanifold. Then $f$ is $C^\infty$ conjugate to an affine Anosov 
   automorphism of a finite factor of a torus.}
\vskip.2cm

In the case of two-dimesional distributions, the additional assumption that 
$\M$ is an infranilmanifold can be replaced by preservation of a volume \cite{F}.

This theorem implies, in particular,  that the conjugacy of such a diffeomorphism 
$f$ to a perturbation $g$ is smooth
if and only if $g$ is also uniformly quasiconformal \cite[Corollary 1.1]{KaS}. 
Hence to establish local rigidity of $f$ it suffices to show that any $C^1$ small 
perturbation with the same periodic data is also uniformly quasiconformal. 
Establishing uniform quasiconformality is also the first step in \cite{L1, L2} 
where smoothness of the conjugacy is obtained directly.

The most general  local rigidity result so far was for uniformly quasiconformal 
Anosov diffeomorphisms $f$ satisfying the following additional assumption: 
\vskip.1cm
\begin{center}
  for any periodic point $p,\;$ 
  $df^m|_{E^s(p)}=a^s(p)\cdot \Id\;$  and 
  $\;df^m|_{E^u(p)}=a^u(p)\cdot \Id$,  
\vskip.1cm
\end{center}
where $m$ is the period of $p$, $E^s$ and $E^u$ are stable and unstable distributions, 
 and $a^s(p)$, $a^s(p)$ are real numbers \cite{L2, KaS}. When one
considers a perturbation $g$ of $f$ with the same periodic data, the derivatives
$dg^m|_{E^s(p)}$ are again multiples of identity. Such a map
preserves any conformal structure, i.e. induces the identity map on the space of conformal
structures on $E^s(p)$.
There is a major difference between this special case and the general one. Indeed, 
for a uniformly quasiconformal Anosov diffeomorphism $f$ these derivatives are conjugate 
to multiples of isometries, and are not necessarily multiples of identity.  This gives 
surprisingly little information about $dg^m|_{E^s(p)}$. It is still conjugate 
to a multiple of isometry, but one has no information on how such conjugacy varies
with $p$. Such a map preserves {\em some} conformal structure, i.e. the induced map
on the space of conformal structures on $E^s(p)$ has a fixed point.  However, it may
significantly affect other conformal structures. In particular, there is no control over
quasiconformal distortion. This makes it difficult to show uniform quasiconformality
of $g$ in the general case. Indeed, so far there have been no results related to
quasiconformality or local rigidity for systems with such periodic data. 
\vskip.1cm
In this paper we introduce some new techniques to study the case of general
quasiconformal periodic data. We establish local, as well as global, rigidity
for such volume preserving systems with two-dimensional stable and unstable 
distributions.

We begin with the study of  transitive Anosov diffeomorphisms for which 
every periodic orbit has only one positive and one negative Lyapunov exponent.  
We note that this assumption does not exclude the possibility of Jordan blocks 
and hence such systems are not necessarily uniformly quasiconformal. Still 
we obtain strong results for these systems which form the basis for further analysis.
In particular, we establish continuity and $C^1$ smoothness of measurable 
invariant conformal structures 
and distributions. We apply these results to volume preserving diffeomorphisms 
with $\dim E^u = \dim E^s =2$ and diagonalizable derivatives of the return 
maps at periodic points.  In addition, we use the Amenable Reduction Theorem.
We thank A. Katok for bringing this result to our attention. In our context this theorem
implies the existence of a measurable invariant conformal structure or a measurable 
invariant one-dimensional distribution for $E^u$ ($E^s$). 
This allows us to establish uniform quasiconformality and hence local and global 
rigidity. We also obtain a local rigidity result for Anosov diffeomorphisms in dimension 3.

\vskip.1cm

We formulate our main results in the next section and prove them in 
Section \ref{proofs}. In Section \ref{preliminaries} we introduce the 
notions used throughout this paper.


\section{Statements of Results} \label{results}

First we consider transitive Anosov diffeomorphisms for which 
every periodic orbit has only one positive and one negative Lyapunov exponent,
i.e. all (un)stable eigenvalues of the derivative of the return map have the same 
modulus.
We do not assume that the derivatives of the return maps at periodic points are
diagonalizable. This class includes systems which are not uniformly quasiconformal
and not locally rigid. However, the following theorem shows that they exhibit 
a variety of useful properties.

\begin{theorem} \label{main} Let $f$ be a transitive Anosov diffeomorphism of 
a compact  manifold $\M$. Suppose that for each periodic point $p$
there is only one positive Lyapunov exponent $\lambda^{(p)}_+$
and only one negative 
Lyapunov exponent $\lambda^{(p)}_-$. Then
\vskip.1cm

\begin{enumerate}

\item Any ergodic  invariant measure for $f$ has only one positive and
only one negative Lyapunov exponent. 
Moreover, for any $\e>0$ there exists $C_{\e}$ such that 
for all $x$ in  $\M$ and  $n$ in $\Z$ 
\vskip.2cm
$K^u(x,n)=${\large$\frac{\max\,\{\,\|\,df^n(v)\,\|\, :\; v\in E^u(x), \;\|v\|=1\,\}}
            {\,\min\,\{\,\|\,df^n(v)\,\|\, :\; v\in E^u(x), \;\|v\|=1\,\}}$}
            $\le C_{\e} e^{\e |n|}$ \; and
\vskip.2cm
$K^s(x,n)=${\large$\frac{\max\,\{\,\|\,df^n(v)\,\|\, :\; v\in E^s(x), \;\|v\|=1\,\}}
            {\,\min\,\{\,\|\,df^n(v)\,\|\, :\; v\in E^s(x), \;\|v\|=1\,\}}$}
            $\le C_{\e} e^{\e |n|}$    \vskip.2cm

\item The stable and unstable distributions $E^s$ and $E^u$ 
        are $C^{1+\beta}$ for some $\beta>0$.
     \vskip.1cm

\item 
There exist $C>0$, $\beta>0$, and $\delta_0>0$ such that for any 
$\delta<\delta_0$, $x,y\in \M$ and $n\in \N$ with $
  \dist\,(f^i(x), f^i(y)) \le \delta$ for $0\leq i \leq n,$
we have 
$$
 \|(df^n_x)^{-1} \circ df^n_y - \Id \,\| \leq C\delta^{\beta}.
$$
     
 \item Any $f$-invariant measurable conformal structure on $E^u$ ($E^s$)
 is $C^1$. Measurability here can be understood with respect 
 to the measure of maximal entropy or with respect to the invariant volume, 
 if it exists.
 \vskip.1cm
       
\item If in addition $f$ is volume-preserving, then any measurable $f$-invariant 
distribution in $E^u$ ($E^s$) defined almost everywhere with respect to the 
volume is $C^1$.

\end{enumerate}
\end{theorem}

\begin{remark} \label{identify}
To consider the composition of the derivatives in (3) we identify the tangent 
spaces at nearby points $x$ and $y$. This can be done by fixing a smooth background
Riemannian metric and using the parallel transport along the unique shortest 
geodesic connecting $x$ and $y$. Such identification can be adjusted using
projections to preserve the Anosov splitting (or an invariant distribution). In this
case the identification will be as regular as the Anosov splitting (or the invariant 
distribution).
\end{remark}

Next we consider the case when $\dim E^u = \dim E^s =2$. 
Now we assume that the matrix of the derivative the return map 
is diagonalizable. We note that the example given by de la Llave 
in \cite{L2} shows that this assumption is necessary in the following 
theorem and its corollary.  The  theorem establishes global rigidity 
for systems with conformal periodic data, and the corollary yields local rigidity.

\begin{theorem} \label{dim4} Let $\M$ be a compact  4-dimensional manifold 
and let $f:\M\to\M$ be a transitive $C^\infty$ Anosov
diffeomorphism with 2-dimensional stable and unstable distributions.
Suppose that for each periodic point $p$ the derivative of the return map 
is diagonalizable over $\C$ and its eigenvalues $\lambda_1^{(p)}, \lambda_2^{(p)},
\lambda_3^{(p)}, \lambda_4^{(p)}$ satisfy 
$$
|\lambda_1^{(p)}|=|\lambda_2^{(p)}|,\quad |\lambda_3^{(p)}|=|\lambda_4^{(p)}|, 
\;\text{ and }\;\; \,|\lambda_1^{(p)}  \lambda_2^{(p)} \lambda_3^{(p)}\lambda_4^{(p)}|=1.
$$
Then $f$ is uniformly quasiconformal and its finite cover is $C^\infty$ 
conjugate to an Anosov automorphism of $\,\T^4$.

\end{theorem}

Note that  the last condition on the eigenvalues is equivalent to 
the fact that  $f$ preserves a smooth volume \cite[Theorem 19.2.7]{KH}.

\begin{corollary} \label{perturbation} 
Let $\M$ and $f$ be as in Theorem \ref{dim4} and let $g : \M \to \M$ be a 
$C^\infty$ Anosov diffeomorphism conjugate to $f$ by a homeomorphism $h$.
Suppose that for any point $p$ such that $f^m(p)=p$, the derivatives $df^m_p$ 
and $dg^m_{h(p)}$ have the same Jordan normal form. 
Then $h$ is a $C^\infty$ diffeomorphism, i.e. $g$ is $C^\infty$ conjugate to $f$. 

\end{corollary}

This corollary applies to any Anosov automorphism of $\,\T^4$ which 
is diagonalizable over $\C$ and whose eigenvalues satisfy 
$\;|\lambda_1|=|\lambda_2| <1<|\lambda_3|=|\lambda_4|$.
\vskip.2cm

We also obtain local rigidity in dimensional three even though there is no 
global rigidity result for this case.

\begin{corollary} \label{dim3} 
Let $f$ be a $C^\infty$ volume-preserving Anosov diffeomorphism of a 
3-dimensional manifold $\M$.  Suppose that for each periodic point $p$ 
the derivative of the return map is diagonalizable over $\C$ and two of 
its eigenvalues have the same modulus.
Let $g$ be a $C^\infty$ diffeomorphism of $\M$ which is  $C^1$ close to $f$ 
and has  the same periodic data. Then $g$ is $C^\infty$ conjugate to $f$. 
\end{corollary}

We note that the matrix of an Anosov automorphism of $\T^3$ has either 
a pair of complex eigenvalues and a real eigenvalue or three real eigenvalues 
of different moduli. Corollary \ref{dim3} applies, in particular, to the former case. 
The local rigidity for the latter case with $C^{1+\text{H\"older}}$ smoothness of 
the conjugacy was recently proved by Gogolev and Guysinsky \cite{GG}. 
Thus Anosov automorphisms of $\T^3$ are locally rigid:

\begin{corollary} \label{dim3auto} 
Let $f$ be an Anosov automorphism of $\,\T^3$ and $g$ be a $C^\infty$ 
diffeomorphism of $\,\T^3$ which is $C^1$ close to $f$ and has the same 
periodic data. Then $g$ is $C^{1+\text{H\"older}}$ conjugate to $f$. 
\end{corollary}


\section{preliminaries} \label{preliminaries}
In this section we briefly introduce the main notions used throughout 
this paper.

\subsection{Anosov  diffeomorphisms}
 
Let $f$ be a diffeomorphism of a compact Riemannian manifold $\M$. 
The diffeomorphism $f$ is called Anosov if there exist a decomposition 
of the tangent bundle $T\M$ into two $f$-invariant 
continuous subbundles $E^s$ and $E^u$, and constants $C>0$, 
$0<\lambda<1$ such that for all $n\in \N$,
  $$
  \| df^n(v) \| \leq C\lambda^n \| v \|
      \;\text{ for }v \in E^s \quad\text{and}\quad
  \| df^{-n}(v) \| \leq C\lambda^n \| v \|
      \;\text{ for }v \in E^u. 
 $$
The distributions $E^s$ and $E^u$ are called stable and unstable. 
It is well-known that these distributions are tangential to the foliations 
$W^s$ and $W^u$ respectively (see, for example \cite{KH}). The leaves of these 
foliations are $C^\infty$ injectively immersed Euclidean spaces, but in general 
the distributions $E^s$ and $E^u$  are only H\"older continuous transversally to 
the corresponding foliations.

\subsection{Uniformly quasiconformal diffeomorphisms} \label{prelim qc}

Let $f$ be an Anosov diffeomorphism of a compact Riemannian manifold 
$\M$. We say that $f$ is uniformly quasiconformal on 
the unstable distribution if the quasiconformal distortion
\vskip.1cm
\begin{center}
$ K^u(x,n)=${\large $\frac{\max\,\{\,\|\,df^n(v)\,\|\, :\; v\in E^u(x), \;\|v\|=1\,\}}
            {\,\min\,\{\,\|\,df^n(v)\,\|\, :\; v\in E^u(x), \;\|v\|=1\,\}}$}
\end{center}
\vskip.1cm
is uniformly bounded for all $n\in\Z $ and $x\in \M$. If $K^u(x,n)=1$ for all $x$ 
and $n$, then $f$ is conformal on $E^u$. Similarly, one can define the corresponding 
notions for $E^s$, or any other continuous invariant distribution. If a diffeomorphism 
is uniformly quasiconformal (conformal)
on both $E^u$ and $E^s$ then it is called uniformly quasiconformal (conformal). 
An Anosov toral automorphism is uniformly quasiconformal if and only if its matrix
is diagonalizable over $\C$ and all its (un)stable eigenvalues are equal in modulus.

Unlike the notion of conformality, the weaker notion of uniform quasiconformality 
does not depend on the choice of a Riemannian metric on the manifold. However,
any transitive uniformly quasiconformal Anosov diffeomorphism is conformal
with respect to some continuous Riemannian metric \cite[Theorem 1.3]{S}.

\subsection{Conformal structures} \label{conf}

A conformal structure on $\R^n$, $n\geq 2$, is a class of proportional 
inner products. The space $\c^n$ of conformal structures on $\R^n$
identifies with the space of real symmetric positive definite $n\times n$ 
matrices with determinant 1, which is isomorphic to $SL(n,\R) /SO(n,\R)$. 
$GL(n,\R)$ acts transitively on $\c^n$
via
$$
X[C] = (\det X^TX)^{-1/n}\; X^T C \, X, \quad 
\text{ where } \; X\in GL(n,\R), \;\text{ and }\; C \in \c^n.
$$ 
It is known that $\c^n$ becomes a Riemannian symmetric space of non-positive 
curvature when equipped with a certain $GL(n,\R)$-invariant metric. The distance
to the identity in this metric is given by
$$
\dist (\Id , C) = \frac {\sqrt{n}}2 \left( (\log \lambda_1)^2 + 
\dots +(\log \lambda_n)^2 \right) ^{1/2},
$$
where $\lambda_1, \dots, \lambda_n$ are the eigenvalues of $C$.
The distance between two structures $C_1$ and $C_2$ can be computed 
as $\dist (C_1, C_2)= \dist (\Id, X[C_2])$,  where $ X[C_1]=\Id.$ 
We note that on any compact subset of $\c^n$ this distance is bi-Lipschitz 
equivalent to the distance induced by the operator norm on matrices.
\vskip.2cm

Now, let $f$ be a diffeomorphism of a compact manifold $\M$ 
and let $E\subset T\M$ be a subbundle invariant under $df$ with 
$\dim E  \geq 2$.  A conformal structure  on $E(x) \subset T_x\M$ 
is a class of proportional  inner products on $E(x)$.  
Using a background Riemannian metric, we can identify an inner product
with a symmetric linear operator with determinant 1 as before. For each $x\in \M$, 
we denote the space of conformal structures on $E(x)$ by $\c(x)$.
Thus we obtain a bundle $\c$ over $\M$ whose fiber over $x$ is $\c(x)$. 
We equip the fibers of $\c$ with the Riemannian metric defined above.  
A continuous (measurable) section of $\c$ is called a continuous 
(measurable) conformal structure on $E$. A measurable conformal 
structure $\tau$ on $E$ is called bounded
if the distance between $\tau(x)$ and $\tau_0(x)$ is uniformly 
bounded on $\M$ for a continuous conformal structure $\tau_0$ on $E$.    

The diffeomorphism $f$ induces a natural pull-back action $F$ on 
conformal structures as follows. For a conformal structure 
$\tau(fx)\in \c(fx)$, viewed as the linear operator on $E(fx)$, 
$\;F_x(\tau(fx))\in \c(x)$ is given by
$$
  F_x(\tau(fx))= \left( \det \, ((df_x)^* \circ df_x) \right)^{-1/n} 
  (df_x)^* \circ \tau(fx) \circ df_x,
$$
where $(df_x)^*:\; T_{fx}\M \to T_{x} \M$ denotes the conjugate operator of $df_x$. 
We note that $F_x: \c_{fx}\to \c_{x}$ is an isometry between the fibers
$\c(fx)$ and $\c(x)$. 

We say that a conformal structure $\tau$ is $f$-invariant
if $F(\tau) = \tau$.
For an Anosov diffeomorphism $f$, a subbundle $E$ can carry 
an invariant conformal structure only if $E\subset E^s$ or $E\subset E^u$.
Clearly, a diffeomorphism is conformal with respect to a Riemannian 
metric on $E$ if and only if it preserves the conformal structure associated 
with this metric. If $f$ preserves a continuous or bounded conformal
structure on $E$ then $f$ is uniformly quasiconformal on $E$.
If  $f$ is a transitive Anosov diffeomorphism and $E$ is H\"older continuous
then the converse is also true: if $f$ is uniformly quasiconformal on $E$ then 
$f$ preserves a continuous conformal structure on $E$ 
(see Theorem 1.3 in \cite{S} and Theorem 2.7 in \cite{KaS2}).


\section{Proofs}  \label{proofs}

\subsection{Proof of Theorem \ref{main}} $\;$
Parts (1)-(5) of the theorem are proven in Propositions 
4.1-4.5 respectively. The propositions are somewhat more 
detailed than the corresponding statements in the theorem.

We would like to point out that while the statements (4) and (5) of the theorem 
look similar, their proofs are completely different. The derivative of $f$ induces an isometry between
the spaces of conformal structures at $x$ and $f(x)$, but these spaces are 
not compact. On the other hand, the Grassman manifold of the subspaces 
at $x$ is compact, but the induced map between the manifolds at $x$ and $f(x)$ 
is not an isometry. This calls for different approaches. 
Moreover, the continuity of 
a measurable invariant conformal structure holds in greater generality. Its proof 
relies only on
the statement (3) of the theorem. It holds, for example, for {\em any} $C^1$ small
perturbation of a conformal Anosov automorphism, with no assumption on the
coincidence of periodic data \cite[Lemma 5.1]{KaS}. In contrast, the proof of 
continuity of a measurable invariant distribution relies on the statement (1) of the 
theorem which requires the coincidence of periodic data. In fact, one may expect
that a typical small perturbation of a conformal Anosov automorphism has simple
Lyapunov exponents and measurable Lyapunov distributions which are not continuous.


\begin{proposition} \label{e-conf}
Let $f$ be a transitive Anosov diffeomorphism of 
a compact  manifold $\M$. Suppose that for each periodic point
there is only one positive Lyapunov exponent (i.e. all unstable 
eigenvalues of the derivative of the return map have the same 
modulus). 

Then any ergodic  invariant measure for $f$ has only one positive Lyapunov 
exponent. Moreover, for any $\e>0$ there exists $C_{\e}$ such that 
$$
  K^u(x,n) \le C_{\e} e^{\e |n|} \;\text{ for all } x \text{ in } \M \text{ and } n \text{ in }\Z.
  $$
\end{proposition}

A similar statement holds for the stable distribution.

\proof
The fact that any ergodic  invariant measure for $f$ has only one 
positive Lyapunov exponent follows from the result of W. Sun and Z. Wang:
\vskip.2cm

 \cite[Theorem 3.1]{SW} {\it Let $\M$ be a $d$-dimensional Riemannian manifold.
Let $f:\M \to \M$  be a $C^{1+\beta}$ diffeomorphism and let $\mu$ be an 
ergodic hyperbolic measure with Lyapunov exponents $\lambda_1 \le ... \le
\lambda_d.$ Then the Lyapunov exponents
of $\mu$ can be approximated by the Lyapunov exponents of hyperbolic 
periodic orbits, more precisely for any $\e>0$
there exists a hyperbolic periodic point $p$ with Lyapunov exponents
$\lambda_1^{(p)} \le ... \le \lambda_d^{(p)}$ such that 
$|\lambda_i-\lambda_i^{(p)}|<\e$ for $i=1, \dots , d$.}
\vskip.2cm

To establish the estimate for the quasiconformal distortion $K^u(x,n)$,
we use the following result.
\vskip.2cm

 \cite[Proposition 3.4]{RH}
{\it Let $f : \M \to \M$ be a continuous map of a compact metric space.
 Let $a_n : \M \to \R$, $n \geq 0$ be a sequence of continuous functions 
such that 
\begin{equation} \label{3.1}
    a_{n+k} (x) \le a_n (f^k (x)) + a_k (x) 
    \;\text{ for every }x \in \M,\;\; n, k \geq 0
\end{equation}    
 and such that there is a sequence of continuous functions $b_n$, $n \geq 0$
satisfying 
\begin{equation} \label{3.2}
    a_n (x) \le a_n (f^k (x)) + a_k (x) + b_k (f^n (x))
     \;\text{ for every } x \in \M,\;\; n, k \geq 0.  
\end{equation}   
If  $\;\inf _n \left( \frac1n \int _\M a_n d \mu \right) < 0\;$ 
for every ergodic $f$-invariant measure, then there
is $N \geq 0$ such that 
$a_N (x) < 0$ for every $x \in \M$.}
\vskip.2cm

We take $\e >0$ and apply this result to $\;a_n(x)=\log K^u(x,n) - \e n.$ 
It is easy to see that the quasiconformal distortion is submultiplicative, i.e. 
$$
 K^u(x, n+k) \le K^u(x,k) \cdot K^u(f^kx, n)\;\text{ for every }x \in \M,\;\; n, k \geq 0.
 $$
Hence the functions $a_n$ satisfy  \eqref{3.1}. It is straightforward to verify that 
$$
K^u(x, n+k) \geq K^u(n,x) \cdot (K^u(f^nx,k))^{-1}
$$
This inequality implies $\, a_{n+k}(x) \geq a_n(x)-b_k(f^nx)$ where 
$b_n(x)=\log K^u(x,n) +\e n$. Taking into account \eqref{3.1} we obtain
\eqref{3.2}.

Since $a_n$ satisfy  \eqref{3.1}, the Subadditive Ergodic Theorem implies
that for every $f$-invariant ergodic measure $\mu$ 
$$
\lim _{n\to \infty} \,\frac 1n {a_n(x)} =\,  \inf _n \,\frac1n \int _\M a_n d \mu  
\quad \text{ for } \mu \; a.e. \; x \in \M.
$$
Since  $\mu$ has only one positive
Lyapunov exponent, it is easy to see that 
$$
\lim _{n\to \infty} \frac 1n \log {K^u(x,n)} =  0  \;\text{ and hence }
\lim _{n\to \infty} \frac 1n {a_n(x)} = -\e<0
\; \text{ for } \mu \; a.e. \; x \in \M.
$$
Thus all assumptions of the proposition above are satisfied
and hence for any $\e >0$ there exists $N_\e$ such that $a_{N_\e}(x)<0$,
i.e. $K^u(x,N_\e) \le e^{\e N_\e}$ for all $x \in \M$. We conclude that
$ K^u(x,n) \le C_{\e} e^{\e n}$ for all $x$ in $\M$  and $n$ in $\N$,
where $C_\e = \max K^u(x,n)$ with the maximum taken over all $x \in \M$ 
and $1 \le n < N_\e$.  Since $K^u(x,n)= K^u(f^nx,-n)$ we obtain 
$ K^u(x,n) \le C_{\e} e^{\e |n|}$ for all $x$  in $\M$ and $n$ in $\Z$.
$\QED$
\vskip.2cm


\begin{proposition} \label{splitting}
Let $f$ be an Anosov diffeomorphism of 
a compact  manifold $\M$. Suppose that 
 for any $\e>0$ there exists $C_{\e}$ such that 
$$
  K^s(x,n) \le C_{\e} e^{\e n} \;\text{ for all } x \text{ in } \M \text{ and } n \text{ in }\N.
 $$
Then the dependence of $E^u(x)$ on $x$ is $C^{1+\beta}$ for some $\beta>0$.\\ 
Similarly, if $K^u(x,n) \le C_{\e} e^{\e n}$ for all $x\in \M,\, n\in \N,\,$ then $E^s$ is $C^{1+\beta}$.
\end{proposition}


\proof

We use the $C^r$ Section Theorem of M. Hirsch, C. Pugh, and M. Shub 
(see Theorems  3.1, 3.2, 3.5,  
and Remarks 1 and 2 after Theorem 3.8 in  \cite{HPS}).  
\vskip.3cm

\cite[$C^r$ Section Theorem]{HPS} \; {\it  Let $f$ be a $C^r$, $r\geq 1$, 
diffeomorphism of a compact $C^r$ manifold $\,\M$. 
Let $\B$ be a $C^r$  finite-dimensional normed vector bundle 
over $\M$ and let $\B_1$ be the corresponding bundle of unite balls in $\B$.
Suppose that $F:\B \to \B$ (or $F:\B_1 \to \B_1$) is a  $C^r$  
extension of $f$ which contracts fibers, i.e. for any $x \in \M$ 
and any $v,w \in \B(x)$ (resp. $v,w \in \B_1(x)$)
$$
   \|F(v)-F(w)\|_{fx} \le k_x \|v-w\|_x \;\text{ with }\, \sup_{x \in \M}  k_x <1.
$$
Then there  exists a unique continuous $F$-invariant  section of $\B$ 
(resp. $\B_1$). If also
$$
   \sup_{x \in \M}  k_x \a_x ^r <1 \;\text{ where }
    \a_x = \|(df_x)^{-1}\|, \label{C^r eq}
$$
then the unique invariant section is  $C^r$}. 
\vskip.2cm 

Since the stable and unstable distributions are a priori only H\"older continuous,
we take $C^2$ distributions $\bar E^u$ and $\bar E^s$  which are close to 
$E^u$ and $E^s$ respectively.  We consider a vector bundle $\B$ whose fiber 
over $x$ is the set of linear operators from $\bEu (x)$ to  $\bEs (x)$. We endow 
the fibers of $\B$ with the standard operator norm. We fix a sufficiently large $n$ 
and consider the graph transform action $F$ induced by the differential of $f^n$. 
More precisely, if $A : \bEu (x) \to \bEs (x)$ is in $\B(x)$ and $L \subset T_x\M$ 
is the graph of this operator then $F(A)\in \B(f^n x)$ is defined to be the operator
from $\bEu (f^n x)$ to  $\bEs (f^n x)$ whose graph is $df^n_x (L)$. We note that
$F(A)$ is well-defined as long as $df^n_x (L)$ is transversal to $\bEs (f^n x)$. 
Let us denote
  $$  \,\, l_x = \min \,\{\,\|df^n_x (v) \| : \; v \in E^u(x),\;  \|v\|=1\,\} $$
  $$ m_x = \min \,\{\,\|df^n_x (v) \| : \; v \in E^s(x), \; \|v\|=1\,\} $$
  $$ \, M_x = \max \,\{\,\|df^n_x (v) \| : \; v \in E^s(x),\;  \|v\|=1\,\}.$$
In the case when 
$$\bar E^u(x)= E^u(x), \quad \bar E^s (x)=E^s (x), \quad 
    \bar E^u (f^n x)= E^u (f^n x), \quad \bar E^s (f^n x)=E^s (f^n x)$$
map $F$ is defined on the whole fiber $\B(x)$ and $F_x:\B(x)\to \B(f^nx)$ is the linear map 
$$F_x (A) = df^n_x |_{E^s (x)} \circ A \circ (df^n_x |_{E^u (x)})^{-1}, \quad \text{with} \qquad
\|F_x \| \le M_x/l_x <1.$$
In particular, $F_x (\B_1(x)) \subset \B_1(f^n x)$. Clearly, the same inclusion holds if 
$\bar E^u$ and $\bar E^s$ are close enough to $E^u$ and $E^s$ respectively. Thus 
$F:\B_1 \to \B_1$ is a well-defined $C^2$ extension of $f$.
In general, $F_x:\B_1(x)\to \B_1(f^nx)$ is an algebraic map which depends 
continuously on the choice of $\bar E^u$ and $\bar E^s$ at $x$ and $f^n x$. 
Thus by choosing $\bar E^u$ and $\bar E^s$ sufficiently close to $E^u$ and 
$E^s$ one can make $F_x$ sufficiently $C^1$-close to the linear map above.
Therefore, $F_x$ is again a contraction with $k_x \approx M_x/l_x$.

We will now apply the $C^r$ Section Theorem to show the smoothness
of the invariant section for $F$. By uniqueness, the graph of this invariant 
section is the distribution $E^u$. The extension $F$ satisfies
  $$
   k_x \approx M_x/l_x <1 \quad \text{and}  \quad  \a_x = \|(df^n|_x)^{-1}\|=1/m_x     
  $$
and thus
$$
    k_x \a_x ^{1+\beta} \, \approx \,\frac{M_x}{l_x m_x^{1+\beta}} \, \le
    \, \frac{1}{l_x\cdot m_x^{\beta}} 
    \cdot \frac{M_x}{m_x} \,\le\, \frac{K^s(x,n)}{l_x \cdot m_x^{\beta}}
     \,\le\, \frac{C_{\e} e^{\e n}}{C l^n \cdot (\inf_x \{m_x\})^{\beta}}
 $$
for some $C>0$ and $l>1$. 

We can take $\e$ so small that $e^\e/l <1$ and then take $n$ so large
that ${C_{\e} e^{\e n}}/{C l^n} <1$. Then the right hand side
is less than 1 for some $\beta >0$. Once $n$ and $\beta$ are chosen, 
we can take $\bEu$ and $\bEs$ close enough to $E^u$ and $E^s$ to
guarantee that $ \sup_{x \in \M}  k_x \a_x ^{1+\beta} <1$. Hence, by the 
$C^r$ Section Theorem, the distribution $E^u$ is $C^{1+\beta}$.
$\QED$
\vskip.2cm


\begin{proposition} \label{df} 
 Let $f$ be an Anosov diffeomorphism of a compact manifold $\M$  such that for any $\e>0$ there exists $C_{\e}$ such that 
\begin{equation} \label{KuKs}
  K^u(x,n) \le C_{\e} e^{\e n} \;\;\text{ and }\;\;  K^s(x,n) \le C_{\e} e^{\e n}
\end{equation}
 for all  $x$ in  $\M$  and  $n$  in $\N$.

Then there exist $C>0$, $\beta>0$, and $\delta_0>0$ such that for any 
$\delta<\delta_0$, $x,y\in \M$ and $n\in \N$ with $
  \dist\,(f^i(x), f^i(y)) \le \delta$ for $0\leq i \leq n,$
we have 
\begin{equation}  \label{dxdy}
 \|(df^n_x)^{-1} \circ df^n_y - \Id \,\| \leq C\delta^\beta.
\end{equation}
  
\end{proposition}

To consider the composition of the derivatives we identify the tangent spaces 
at nearby points preserving the Anosov splitting as in Remark \ref{identify}. 
Since the Anosov splitting is $C^{1+\beta}$  by Proposition \ref{splitting}, this 
identification is also $C^{1+\beta}$. 

\proof We will use non-stationary linearizations along stable and unstable 
manifolds given by the following proposition.

\vskip.2cm

\cite[Proposition 4.1]{S} 
{\it Let $f$ be a diffeomorphism of a compact Riemannian manifold $\M$,\, 
and\, let\, $W$ be a continuous invariant foliation with $C^\infty$
leaves.\; Suppose that $\| df|_{TW}\|<1$, and there exist 
$C>0$ and $\e>0$ such that for any $x\in \M$ and $n\in\mathbb N$, 
\begin{equation}  \label{2-pinching}
  \|\,(df^n|_{T_xW})^{-1}\,\| \cdot \|\,df^n|_{T_xW}\,\|^2 
  \leq C(1-\e)^n.  
\end{equation}
Then for any $x\in \M$ there exists a $C^\infty$ diffeomorphism 
$h^{_W}_x: W(x) \to T_xW$ such that

(i)  $h^{_W}_{fx}\circ f=df_x \circ h^{_W}_x,$

(ii)  $h^{_W}_x(x)=0$ and $(dh^{_W}_x)_x$ is the identity map,

(iii) $h^{_W}_x$ depends continuously on $x$ in $C^\infty$ topology.}
\vskip.2cm
 
 Clearly, condition \eqref{KuKs} implies \eqref{2-pinching} with $W=W^s$
 and its analogue for the unstable distribution.
 Thus we obtain linearizations $h^s_x:W^s(x) \to E^s(x)$ and 
$h^u_x:W^u(x) \to E^u(x)$.
 Then we construct a local linearization $h_x:U_x \to T_x\M$,  
 where $U_x$ is a small open neighborhood of $x\in \M$ 
 as follows: 
 $$ 
  h_x|_{W^u(x)}=h^u_x, \quad h_x|_{W^s(x)}=h^s_x,  
 $$
and for $y\in W^u(x)\cap U_x$ and $z\in W^s(x)\cap U_x$ we set
 $$
  h_x([y,z]) = h^u_x(y) + h^s_x(z),   
 $$
where $[y,z]=W_{\text{\text{loc}}}^s(y) \cap W_{\text{\text{loc}}}^u(z)$. It is easy to see that $h$
satisfies conditions (i) and (ii). Since the maps $h^s_x$ and $h^u_x$ are
$C^\infty$ and the local product structure is $C^{1+\beta}$, the maps 
$h_x$ are $C^{1+\beta}$ with uniform H\"older constant for all $x$.
\vskip.1cm

Since the differential $df_x$ is the direct sum of the stable differential 
$df|_{E^s(x)}$ and the unstable differential $df|_{E^u(x)}$, it suffices 
to prove the proposition for these restrictions. We will give a proof  for 
the unstable differential, the other case is similar. 

If $\delta_0$ is small enough, there exists a unique
point $z \in W_{\text{\text{loc}}}^u(x)\cap W_{\text{\text{loc}}}^s(y)$ with
 $$
   \dist(f^i(x), f^i(z))<C_1\delta \;\text{ and }\;\dist(f^i(z), f^i(y))<C_1\delta
   \quad\text{ for }\;0\leq i \leq n.
 $$
Thus it is sufficient to prove the proposition for $x$ and $y$ lying on the same 
stable or on the same unstable manifold. First we consider the case when 
$y\in W^s(x)$. When considering the unstable differential this case appears 
more difficult, but there will be little difference in our argument.  

Since $y \in W^s(x)$ and $\dist\,(x, y) \le \delta$ we obtain
$\dist\,(f^n x, f^n y) \le \delta \cdot C_2 \gamma ^n$ for some positive 
$\gamma <1$. Using the linearizations $h_x$ and $h_{f^nx}$ we can write
in a neighborhood of $x$
$$
f^n= (h_{f^nx})^{-1} \circ df^n _x \circ h_x.
$$ 
Differentiating at $y$ we obtain
$$
df^n_y= ((dh_{f^nx})_{f^ny})^{-1} \circ df^n _x \circ (dh_x)_y
$$ 
and restricting to $E^u$
\begin{equation} \label{1}
       df^n|_{E^u(y)}= 
       \left( (dh_{f^nx})|_{E^u(f^ny)}\right) ^{-1} \circ df^n |_{E^u(x)} \circ (dh_x)|_{E^u(y)}.
\end{equation} 
Since the linearizations are $C^{1+\beta}$ we obtain that
$$
(dh_z)_w =\Id + R \quad \text{ with } \;\; \| R \| \le C_3 \, \dist (z,w)^\beta
$$ 
for all $z \in \M$ and all $w$ close to $z$. Therefore we can write
\begin{equation} \label{2}
     (dh_{f^nx})|_{E^u(f^ny)} =
      \Id + R_1 \quad \text{ and } \;\; (dh_x)|_{E^u(y)}=\Id + R_2,
\end{equation} 
where
\begin{equation} \label{3}
  \| R_2 \| \le C_3 \, \dist (x,y)^\beta \le C_3 \, \delta^\beta
\end{equation} 
and
$$
  \| R_1 \| \le C_3 \, \dist (f^nx,f^ny)^\beta 
  \le C_3 ( C_2 \,\delta \gamma ^n)^\beta \le C_4 \,\delta^\beta \gamma ^{\beta n}
$$
We also observe that 
\begin{equation} \label{4}
    \left( (dh_{f^nx}) |_{E^u(f^ny)}  \right) ^{-1} =(\Id + R_1)^{-1} =\Id + R_3
\quad \text{ with } 
\end{equation} 
\begin{equation} \label{5}
    \| R_3 \| \le \frac{\| R_1 \|}{1- \| R_1\|} \le 
    2 C_4 \,\delta^\beta \gamma ^{\beta n}
\end{equation}  
provided that $\| R_1 \| <1/2$.
Combining \eqref{1}, \eqref{2}, and \eqref{4}, we can write
$$
(df^n |_{E^u(x)})^{-1} \circ df^n |_{E^u(y)}= 
(df^n |_{E^u(x)})^{-1} \circ (\Id + R_3) \circ df^n |_{E^u(x)} \circ (\Id + R_2)
$$
$$ 
=\Id + R_2 + (df^n |_{E^u(x)})^{-1} \circ R_3 \circ df^n  |_{E^u(x)} \circ (\Id + R_2)
$$
Therefore, we can estimate
$$
\|(df^n |_{E^u(x)})^{-1} \circ df^n |_{E^u(y)} - \Id \|  \le  
$$
$$
\| R_2 \| + \|(df^n|_{E^u(x)})^{-1}\| \cdot \| R_3 \| \cdot \|  df^n |_{E^u(x)} \| \cdot \| \Id + R_2 \|
$$
We note that $\|(df^n|_{E^u(x)})^{-1}\| \cdot  \|  df^n |_{E^u(x)} \| \le K^u (x,n) \le C_\e e^{\e n}$.
Finally, using estimates  \eqref{3} and \eqref{5} for $\| R_2 \|$ and $\| R_3 \|$, we obtain
$$
\|(df^n |_{E^u(x)})^{-1} \circ df^n |_{E^u(y)} - \Id \|  \le
C_3 \, \delta^\beta + 
2 C_4 \,\delta^\beta \gamma ^{\beta n} \cdot C_\e e^{\e n} \cdot (1+ C_3 \, \delta^\beta)
$$
If $\e$ is chosen sufficiently small so that $|\gamma ^{\beta}  e^{\e}|<1$, then the term
$\gamma ^{\beta n} \cdot C_\e e^{\e n}$ is uniformly bounded in $n$. So we conclude
that
$$
\|(df^n |_{E^u(x)})^{-1} \circ df^n |_{E^u(y)} - \Id \|  \le
C_5 \, \delta^\beta. 
$$
\vskip.1cm

To complete the proof of the proposition it remains to consider the case when 
$y\in W^u(x)$. We use the same notation as in the previous case and indicate
the necessary changes. In this case we have $\dist\,(f^n x, f^n y) \le \delta$ and
$\dist\,(x, y) \le \delta \cdot C_2 \gamma ^n$ for some positive $\gamma <1$.
Therefore, 
\begin{equation} \label{6}
 \| R_3 \| \le \frac{\| R_1 \|}{1- \| R_1\|} \le 2 C_3 \, \delta^\beta
\; \text{ and } \;
 \| R_2 \| \le C_3 \, \dist (x,y)^\beta \le C_6 \, \delta^\beta \gamma ^{\beta n}.
\end{equation} 
Now we can write
$$
df^n |_{E^u(y)}  \circ (df^n |_{E^u(x)})^{-1}= 
(\Id + R_3) \circ df^n |_{E^u(x)} \circ (\Id + R_2) \circ (df^n |_{E^u(x)})^{-1}
$$
$$ 
=\Id + R_3 +  (\Id + R_3) \circ (df^n |_{E^u(x)})^{-1} \circ R_2 \circ df^n  |_{E^u(x)}.
$$
Finally, using the new estimates \eqref{6} for $\| R_2 \|$ and $\| R_3 \|$, we obtain
similarly to the previous case that
$$
\| df^n |_{E^u(y)}  \circ (df^n |_{E^u(x)})^{-1} - \Id \|  \le
C_7 \, \delta^\beta. 
$$
This estimate clearly implies a similar H\"older estimate for 
$\|(df^n |_{E^u(x)})^{-1} \circ df^n |_{E^u(y)} -\Id \|\,$ 
and concludes the proof of the proposition.
$\QED$
\vskip.2cm


\begin{proposition} \label{structure}$\;$
  
(i) Let $f$ be a a transitive Anosov diffeomorphism of a compact manifold $\M$
and let $E$ be a H\"older continuous invariant distribution.
Suppose that  there exist $k>0$, $\delta_0>0$, and $\beta>0$ such that for any 
$\delta<\delta_0$, $x,y\in \M$ and $n\in \N$ such that
$\;\dist\,(f^i(x), f^i(y))<\delta$ for $0\leq i \leq n,\;$ 
$$
  \|\left( df^n|_{E(x)}\right)^{-1} \circ \left( df^n|_{E(y)} \right) - \Id \;\| 
  \leq k\,\delta^\beta.
$$ 
Then any $f$-invariant measurable  conformal structure 
on  $E$ is H\"older continuous. 
\vskip.1cm

(ii) If $f$ is as in Theorem \ref{main}, then 
any $f$-invariant measurable conformal structure on $E^u$ ($E^s$) is $C^1$.

\end{proposition}

{\bf Remark.} As follows from the proof, the measurability of the conformal 
structure in this proposition can be understood with respect to any ergodic 
measure $\mu$ with full support for which the local stable (or unstable) 
holonomies are absolutely continuous with respect to the conditional 
measures on local unstable (or stable) manifolds. Note that invariant 
volume, if it exists, satisfies these properties, and so does the Bowen-Margulis 
measure of maximal entropy, which always exists for a transitive Anosov 
diffeomorphism $f$. More generally, these properties hold for any equilibrium 
(Gibbs) measure corresponding to a H\"older continuous potential.

\proof
{\bf (i)} In this proposition we identify spaces $E(x)$ and $E(y)$ at nearby points 
$x$ and $y$ as in Remark \ref{identify}. Since the distribution $E$ is H\"older 
continuous,  the identification is also H\"older continuous, and hence 
$df|_E$ is  H\"older continuous with respect to the identification. 
The identification allows us to conveniently compare differentials and conformal structures at different points.

For $x\in \M$, we denote by $\tt(x)$ the conformal structure on $E(x)$.
First we estimate the distance  between the 
conformal structures at $x$ and at a nearby point $y\in W^s(x)$.
We use the distance described in Section \ref{conf}.
Let $x_n=f^n(x)$, $y_n=f^n(y)$, and let $F^n_x$ be the isometry 
from $\c(f^nx)$ to $\c(x)$ induced by $df^n|_{E(x)}$.
Since the conformal structure $\tt$ is invariant, 
$\tt(x)=F_x^n (\tt(x_n))$  and $\tt(y)=F_y^n (\tt(y_n))$. 
Using this and the fact that  $F^n$ is an isometry, we obtain
 $$
  \begin{aligned}
  \dist(\tt(x), \tt(y))&=\dist \left( F_x^n (\tt(x_n), F_y^n (\tt(y_n)) \right) \\
  &\leq \dist \left( F_x^n(\tt(x_n), F_y^n(\tt(x_n))  \right) +
         \dist   \left( F_y^n(\tt(x_n), F_y^n(\tt(y_n))  \right)\\
  &=\dist \left( \tt(x_n), ((F_x^n)^{-1} \circ F_y^n)(\tt(x_n))  \right) +
         \dist  \left( \tt(x_n),\tt(y_n)) \right).
 \end{aligned}
 $$

Let $\mu$ be an invariant measure as in the proposition or the remark after it.
Since the conformal structure $\tt$ is measurable, by Lusin's theorem we can 
take a compact set $S\subset \M$ with $\mu(S)>1/2$ on which $\tt$ is uniformly 
continuous and bounded. 

First we show that for $x_n \in S$ the term 
$\dist \left( \tt(x_n), ((F_x^n)^{-1} \circ F_y^n)(\tt(x_n))  \right)$ 
is H\"older in $\dist (x,y)$. For this we observe that the map 
$(F_x^n)^{-1} \circ F_y^n$ is induced by $(df^n_x)^{-1}\circ df_y^n,\,$
and $\,\|(df^n_x)^{-1} \circ df^n_y - \Id \,\| \leq k\cdot \dist(x,y)^\beta$. 
Let $A$ be the matrix of $(df^n_x)^{-1} \circ df^n_y$. Then 
$$
A=\Id+R, \;\text{ where }\,\|R\| \le  k\cdot \dist(x,y)^\beta.
$$ 
Let $C$ be the matrix corresponding to the conformal structure $\tt(x_n)$.
Recall that $C$ is symmetric positive definite with determinant 1.
Thus there exists an orthogonal matrix $Q$ such that 
$Q^T C Q$ is a diagonal matrix whose diagonal entries are the  eigenvalues
$\lambda_i>0$ of $C$. Let $X$ be the product of $Q$ and the diagonal matrix 
with entries $1/\sqrt{\lambda_i}$. Then $X$ has determinant 1
and  $X[C]=X^T C X=\Id$.  Now we estimate
$$
\begin{aligned}
  &\dist \left( \tt(x_n), ((F_x^n)^{-1} \circ F_y^n)(\tt(x_n)) \right) =
               \dist\, (C, A[C])=\;\dist\, (\Id, X[A[C]]) \\
  &=\dist \left( \Id, X^TA^TCAX \right)
 =\dist \left( \Id, X^T(\Id+R^T)C(\Id+R)X \right) \\&=\dist\, (\Id, \Id+B),
 \;\text{ where }\; B=X^TCRX+X^TR^TCX+X^TR^TCRX.
  \end{aligned}
$$
We observe that $\|B\| \le 3 \|X\|^2\cdot \|C\| \cdot \|R\|$ and $\|X\|^2 \le \|C^{-1}\|$,
as follows from the construction of $X$.
Since the conformal structure $\tt$ is bounded on $S$, so are $\|C^{-1}\|$ and $\|C\|$, 
and hence 
$$
\begin{aligned}
  & \|B\| \le 3 \|C^{-1}\|\cdot \|C\| \cdot \|R\| 
  &\le  3 k_1  \|R\| \,\le \, 3 k_1 k \cdot \dist(x,y)^\beta=
  k_2\cdot \dist(x,y)^\beta
 \end{aligned}
$$
In particular, $\|B\|$ is small if $\dist(x,y)$ is sufficiently small. Finally, for $x_n \in S$
we obtain
$$
\begin{aligned}
  & \dist \, ( \tt(x_n),  ((F_x^n)^{-1} \circ F_y^n)(\tt(x_n)) ) =\dist\, (\Id, \Id+B)  \le k_3 \|B\| \le
  k_4\cdot \dist(x,y)^\beta
 \end{aligned}
$$
where constant $k_4$ depend on the set $S$.
We conclude that if $x_n$ is in $S$ then
 $$
 \dist(\tt(x), \tt(y))
  \leq \dist(\tt(x_n),\tt(y_n)))+k_4 \cdot \dist(x,y)^\beta.
 $$ 
 
 It follows from the Birkhoff ergodic theorem that the set of points for which 
the frequency of visiting $S$ equals $\mu(S)>1/2$ has full measure. 
We denote this set by $G$. If both $x$ and $y$ are in $G$, 
then there exists a sequence $(n_i)$ such that $x_{n_i}\in S$ and 
$y_{n_i}\in S$.
If in addition $x$ and $y$ lie on the same stable leaf,
then  $\dist(x_{n_i}, y_{n_i})\to 0$ 
and hence $\dist(\tt(x_{n_i}), \tt(y_{n_i}))\to 0$ by continuity of $\tau$ on $S$. 
Thus, we obtain
 $$
   \dist(\tt(x), \tt(y))\leq  k^s \cdot\dist(x,y)^\beta.
 $$
By a similar argument, 
$\dist(\tt(x), \tt(z))\leq  k^u \cdot\dist(x,z)^\beta$ for any 
two nearby points $x, z\in G$ lying on the same unstable 
leaf. 

Consider a small open set in $\M$ with a product structure.
For $\mu$ almost all local stable leaves, the set of  points of 
$G$ on the leaf has full conditional measure.  
Consider points $x,y \in G$ lying on two such  stable leaves.
Let $H_{x,y}$ be the unstable holonomy map between $W^s(x)$
and $W^s(y)$. Since the holonomy maps are absolutely continuous
with respect to the conditional measures, there exists a point 
$z\in W^s(x)\cap G$ close to $x$ such that $H_{x,y}(z)$ is also in $G$.
By the above argument,  
  $$
   \begin{aligned}
    \dist (\tt(x), \tt(z))&\leq  k^s \cdot\dist(x,z)^\beta,  \\
    \dist(\tt(z), \tt(H_{x,y}(z))) &\leq  k^u \cdot
       \dist(z,H_{x,y}(z))^\beta, \quad\text{and} \\
    \dist(\tau(H_{x,y}(z)),\tau(y))&\leq k^s\cdot
        \dist(H_{x,y}(z),y)^\beta.
   \end{aligned}
  $$
Since the points $x$, $y$, and $z$ are close, it is clear from
the local product structure that
 $$ 
   \dist(x,z)^\beta + \dist(z,H_{x,y}(z))^\beta +
    \dist(H_{x,y}(z),y)^\beta \leq  k_5\cdot\dist(x,y)^\beta.
 $$  
Hence, we obtain $\dist(\tt(x),\tt(y))\leq  k_6
\cdot \dist(x,y)^\beta$ for  all $x$ and $y$ in a set of full measure 
$\tilde G \subset G$.

We can assume that $\tilde G$ is invariant by considering 
$\bigcap_{n=-\infty}^{\infty} f^n(\tilde G)$. Since $\mu$ has
full support the set $\tilde G$ is dense in $\M$. Hence we can extend $\tt$ 
from $\tilde G$ and obtain an invariant H\"older continuous conformal 
structure $ \tau$ on $\M$. 
This completes the proof of the first part of the proposition.
\vskip.2cm

{\bf (ii)} Let $\tau$ be a measurable conformal structure on $E^u$. By (i) 
it is H\"older continuous, and we can show that it is actually smooth as follows.
First we note that by Lemma 3.1 in \cite{S} the conformal
structure $\tau$ is $C^\infty$ along the leaves of $W^u$. The lemma shows 
that the $C^\infty$ linearization $h_x:W^u(x) \to E^u(x)$ as in 
Proposition \ref{df} maps  $\tau$ on $W^u(x)$ to a constant conformal 
structure on $E^u$.

By Proposition \ref{splitting}, the distributions $E^u$ and $E^s$
are $C^{1+\beta}$.
Theorem 1.4 in \cite{S} implies that $\tau$ is preserved by the stable holonomies,
which are $C^{1+\beta}$ (in fact they are $C^\infty$ for $\dim E^u \ge 2$).
Thus $\tau$ is smooth along the leaves of $W^s$.
Now it follows easily, for example from Journe Lemma \cite{J}, that $\tau$ 
is at least $C^1$ on $\M$.
$\QED$
\vskip.2cm


\begin{proposition} \label{distribution} Let $f$ be a volume-preserving 
Anosov diffeomorphism of a 
compact manifold $\M$. Suppose that for each periodic point $p$
there is only one positive Lyapunov exponent $\lambda^{(p)}_+$
and only one negative Lyapunov exponent $\lambda^{(p)}_-$. 
Then any measurable  measurable $f$-invariant distribution in $E^u$ ($E^s$) 
defined almost everywhere with respect to the volume is $C^1$.
\end{proposition}

\proof
Let $E$ be a $k$-dimensional measurable invariant distribution in $E^u$.
We consider fiber bundle $\bM$ over $\M$ whose fiber over $x$ is the Grassman 
manifold $G_x$ of all $k$-dimensional subspaces in $E^u (x)$. The differential 
$df_x$ induces a natural map $F_x : G_x \to G_{fx}$ and we obtain an extension
$\f :  \bM \to \bar \M $ of our system $f:\M\to\M$ given by
$ \f(x,V)=( f(x), F_x(V))$ where  $V\in G_x$. Thus we have the following
commutative diagrams where  $\pi$ is the projection.

\begin{center}
\begin{tabular}{ccc}
$\bM$ & $ \overset{\f}{\longrightarrow}$ & $\bM$ \\
$\;\;\;\downarrow$ \scriptsize{$\pi$} & & $\;\;\;\downarrow$ \scriptsize{$\pi$} \\
$\M$ & $\overset{f}{\longrightarrow}$ & $\M$
\end{tabular} \hskip2cm
\begin{tabular}{ccc}
$T\bM$ & $\overset{d\f}{\longrightarrow}$ & $T\bM$ \\
$\;\;\;\;\downarrow$ \scriptsize{$d\pi$} & & $\;\;\;\;\downarrow$ \scriptsize{$d\pi$} \\
$T\M$ & $\overset{df}{\longrightarrow}$ & $T\M$
\end{tabular}
\end{center}
We note that since the distribution $E^u$ is $C^{1+\beta}$ by Proposition \ref{splitting}, 
the diffeomorphism  $\bar f$ is also $C^{1+\beta}$.
\vskip.1cm

Let $\bEc$ be the distribution in $T\bM$ tangent to the fibers, 
i.e. $\bEc = \ker d\pi.$
Clearly, $\bEc$ is invariant under $d\f$.  

\begin{lemma} \label{neutral_lemma}
There exists $C>0$ such that for any $\bar x \in \bM$,  $n \in \Z$, and
$\e>0$
\begin{equation} \label{neutral}
 \| d\f^n |_{\bEc( \bar x)} \| \le 
C \cdot K^u(x,n) \le C\cdot C_\e e^{\e |n|} , \quad \text{where }\;x=\pi (\bar x).
\end{equation}
\end{lemma}

\proof
Let $w$ and $v$ be two unit vectors in $E^u (x)$.
We denote $D=d f^n_{ x}$. Using the formula
$\; 2<Dw, Dv>=\|Dw\|^2+\|Dv\|^2-\|Dw-Dv\|^2$ for the inner product, we obtain 
$$
(2 \, \sin ({\angle (Dw, Dv)}/2)\,)^2 =2 (1-\cos \angle (Dw, Dv) )=
$$
$$
  2-\frac{2 <Dw, Dv>}{\|Dw\|\cdot \|Dv\|}=
  \frac{2 \|Dw\|\cdot \|Dv\| - \|Dw\|^2 - \|Dv\|^2 + \|Dw-Dv\|^2} {\|Dw\|\cdot \|Dv\|} =
$$  
$$
  \frac{ \|Dw-Dv\|^2 -(\|Dw\| - \|Dv\|)^2 } { \|Dw\|\cdot \|Dv\|} \le
  \frac{ \|D\|^2 \cdot \|w-v\|^2} { \|Dw\|\cdot \|Dv\|} \le K^u(x,n) ^2  \cdot \|w-v\|^2. 
$$ 
Suppose that the angle between the unit vectors $w$ and $v$ is sufficiently small
so that it remains small when multiplied by $K^u(x,n)$. Then we obtain that
the angle $\angle (Dw, Dv)$ is also small and
$$ \angle (Dw, Dv) \le 2 K^u(x,n) \cdot \angle (w, v)$$
Therefore, for any subspaces $V,W \in G_x$ we have
$$ \dist (F^n_x (V), F^n_x (W)) \le 2 K^u(x,n) \cdot \dist (V, W)$$
where the distance between two subspaces is the maximal angle. This
means that $\f^n$ expands the distance between any two nearby points
in the fiber $G_x$ at most by a factor of $2K^u(x,n)$. We note that the 
maximal angle distance is Lipschitz equivalent to any smooth Riemannian 
distance on the Grassman manifold. It follows that there exists $C>0$
such that for any $\bar x \in \bM$ and any $n >0$ 
$$
\| d\f^n |_{\bEc( \bar x)} \| \le C \cdot  K^u(\pi (x),n) \le C\cdot C_\e e^{\e |n|},
$$
where the last inequality is given by Proposition \ref{e-conf}.
The case of $n<0$ can be considered similarly.
$\QED$
\vskip.2cm

We define distributions $\bar E^{uc}= d\pi^{-1} (E^u)$ and 
$\bar E^{sc}= d\pi^{-1} (E^s)$. It follows from the commutative diagram 
that $\bar E^{uc}$ and $\bar E^{sc}$ are  invariant under $d\f$.
Note that $\bEc =\bar E^{uc} \cap \bar E^{sc}$.

\begin{lemma} \label{partially}

 $\f$ is a partially hyperbolic diffeomorphism.  More precisely, there exists a 
 continuous $d\f$-invariant splitting
\begin{equation}\label{part_hyp}
T \bM = \bEu \oplus \bEc \oplus \bEs, 
\end{equation}
where the unstable and stable distributions $\bEu$ and $\bEs$ are contained 
in $\bar E^{uc}$ and $\bar E^{sc}$ respectively and  are expanded/contracted  
uniformly by $\f$. The expansion/contraction in $\bEc$ is uniformly slower than 
expansion in $\bar E^u$ and contraction in $\bar E^s$.

\end{lemma}

\proof We fix a continuous distribution $H$ such that $\bar E^{uc} = \bEc \oplus H$.
We can identify $H$ with $E^u$ via $d\pi$ and introduce the metric on $\bar E^{uc}$
so that the sum is orthogonal. We consider fiber bundle $\B$ whose fiber over 
$x \in \bM$ is the space of linear operators $L_x : H(x) \to \bEc (x)$ equipped
with the operator norm. The differential $d\f$ induces via the graph transform 
the natural extension $F_x: \B (x) \to \B(\f x)$ of $\f$. It is easy to see that $F_x$ 
is an affine map. Lemma \ref{neutral_lemma} implies that the possible expansion 
in $\bEc$ is slower than the uniform expansion in $H \cong E^u$. Similarly to the
proof of Proposition \ref{splitting}, this implies that $F$ is a uniform fiber contraction.
Hence it has a continuous invariant section whose graph gives a continuous
$d\f$ invariant distribution which we denote $\bEu$. Since $d\pi$ induces an 
isomorphism between $\bEu$ and $E^u$ which conjugates the actions of $d\f$
and $df$, we conclude that $d\f$ uniformly expands $\bEu$. The distribution $\bEs$
is obtained similarly.
$\QED$
\vskip.2cm

The theory of partially hyperbolic diffeomorphisms gives the existence of the
unstable and stable foliations $\bWu$ and $\bWs$ with $C^1$ smooth leaves
tangential to the distributions $\bEu$ and $\bEs$ respectively. 
The leaf $\bWu (\bar x)$ is contained in the corresponding center-unstable leaf 
$\bar W^{uc}(\bar x)= \pi ^{-1} (W^{u}(x))$, where $\pi (\bar x) =x$, and projects diffeomorphically to $W^{u}(x)$.  Similarly, $\bWs (\bar x)$ is contained in 
$\bar W^{sc}(\bar x)= \pi ^{-1} (W^{s}(x))$ and projects diffeomorphically to 
$W^{s}(x)$. 

\vskip.3cm

The measurable $f$-invariant distribution $E$ gives rise to the measurable 
$\f$-invariant section $\phi:\M \to \bM$, $\;\phi(x)=(x,E(x))$. 
We will show that $\phi$ is $C^1$, and hence so is $E$.
Let $\mu$ be the smooth $f$-invariant volume on $\M$.
We denote by $\bm$ the lift of $\mu$ to the graph $\Phi$ of $\phi$, i.e.
for a set $X\subset \bM$, $\;\bm(X)=\mu(\pi(X\cap \Phi))$.
Since $\mu$ is $f$-invariant and $\Phi$ is $\f$-invariant,
the measure $\bm$ is also $\f$-invariant.

Now we consider Lyapunov exponents of $\bm$. Lemma \ref{neutral_lemma} 
implies that the Lyapunov exponent of any vector in $\bEc$ is zero. The uniform
contraction in $\bEs$ and expansion in $\bEu$ now imply that the partially 
hyperbolic splitting \eqref{part_hyp} coincides with the Lyapunov splitting 
into the unstable, neutral, and stable distributions for $\bm$.

We recall that by the Ruelle inequality the entropy of any ergodic measure 
invariant under a diffeomorphism is not greater than the sum of its 
positive Lyapunov exponents counted with multiplicities. It is well known that the inequality is, in fact, equality  for an invariant volume. Thus 
\begin{equation} \label{h1}
\begin{aligned}
h_{\bm}(\f) \;&\le \text{ sum of positive exponents of }\bm ,\\
h_{\mu}(f)\; &= \text{ sum of positive exponents of }\mu. 
\end{aligned}
\end{equation}
Lemma \ref{neutral_lemma} implies that 
$\| d\f^n |_{\bar E^{uc} (\bar x)}\| \le \| df^n |_{ E^{u} (x)}\|$ and hence
the largest Lyapunov exponent of $\bm$ is not greater than that of $\mu$.
By Proposition \ref{e-conf} $\mu$ has only one positive exponent,
and since $\dim E^u = \dim \bEu$ we obtain 
\begin{equation} \label{h2}
 \text{sum of positive exponents of }\bm\;  \le
 \text{ sum of positive exponents of }\mu.
\end{equation}
Since $(f,\mu)$ is a factor of $(\f,\bm)$ we have $h_\mu (f) \le h_{\bm}(\f)$.
Combining this with \eqref{h1}  and \eqref{h2} we conclude that
$$
\begin{aligned}
h_{\bm}(\f) \; \le &  \text{ sum of positive exponents of }\bm  \\
 \le &  \text{ sum of positive exponents of }\mu \; = \; h_\mu (f)\; \le\; h_{\bm}(\f),
\end{aligned}
$$
and hence all the inequalities above are, actually, equalities.
In particular, $h_{\bm}(\f)$ equals the sum of the positive exponents of $\bm$.
This implies \cite[Theorem A]{LY1} that the conditional measures on the unstable 
manifolds of $\f$ are absolutely continuous with positive densities. We note that
the $C^2$ smoothness assumption for the diffeomorphism in that theorem is 
only used to establish
that the holonomies of the unstable foliation within a center-unstable set are
Lipschitz. Since $\f$ is $C^{1+\beta}$ partially hyperbolic and dynamically coherent
diffeomorphism satisfying Lemma \ref{neutral_lemma}, 
these holonomies are in fact $C^{1}$ \cite[Theorem 0.2]{BW}. Similarly, 
we establish that the conditional measures on the stable manifolds of $\f$ are 
absolutely continuous with positive densities. We show below 
that the local product structure given by foliations $W^u$ and $W^s$ on $\M$ lifts 
to the local product structure on the graph $\Phi$ in $\bM$ given by foliations 
$\bWu$ and $\bWs$, which easily implies the proposition. We note that the 
absolute continuity of the conditional measures seems essential for this argument.

For a point $x$ in $\M$ we will denote by $\bar x$ the unique point on the graph
$\Phi$ with $\pi (\bar x) =x$. The absolute continuity of the conditional measures
implies that for a $\bm$-typical point $\bar x$ almost every point $\bar y$ with 
respect to a volume on the local stable leaf $\bWs  _{\text{loc}} (\bar x)$ is also typical 
with respect to $\bm$. The projection of the set of such points $\bar y$ 
has full measure in the local stable leaf $W^s _{\text{loc}}  (x)$ with respect to 
a volume on this leaf, which is equivalent to the conditional measure of $\mu$.  
This means that for a $\mu$-typical point $x$ the leaf $\bWs  _{\text{loc}} (\bar x)$
is essentially contained in the graph $\Phi$. More precisely, there is a set
$X$ in $\M$ with $\mu (X)=1$ such that if $x,z \in X$ and $z \in W^s _{\text{loc}}  (x)$
then $\bar z \in \bWs _{\text{loc}}  (\bar x)$. We can also assume that $X$ has a similar
property for the unstable leaves.

Consider nearby points $x,y \in X$. Absolute continuity of the unstable holonomies
for $f$ gives the existence of points $z_1,z_2 \in X$ such
that $z_1 \in W^s _{\text{loc}}  (x)$, $z_2 \in W^s _{\text{loc}}  (y)$, and $z_2 \in W^u _{\text{loc}}  (z_1)$.
By the property of $X$ the corresponding lifts satisfy
$\bar z_1 \in \bar W^s _{\text{loc}}  (\bar x)$, $\bar z_2 \in \bar W^s _{\text{loc}}  (\bar y)$, 
and $\bar z_2 \in \bar W^u _{\text{loc}}  (\bar z_1)$, i.e. $\bar x$ and $\bar y$ can be
connected by local stable/unstable manifolds. Since these manifolds are $C^1$
graphs over the corresponding manifolds in $\M$, we conclude that
$\dist ( \bar x, \bar y) \le C \cdot \dist (x,y)$. This implies that the graph $\Phi$ and 
the section $\phi$ are Lipschitz continuous (up to a change on a set of measure 
zero). Moreover, the graph of $\phi$ restricted to $W^s(x)$ is the manifold 
$\bWs (\bar x$), and hence the restriction is $C^1$. Similarly, the restriction of $\phi$ 
to $W^u(x)$ is $C^1$ for any $x \in \M$. It follows, for example from 
Journe Lemma \cite{J}, that $\phi$ is $C^1$ on $\M$.
$\QED$
\vskip.2cm 
This completes the proof of Theorem \ref{main}.


\subsection{Proof of Theorem \ref{dim4}}$\;$ 
We will use the following particular case of Zimmer's Amenable 
Reduction Theorem (see \cite[Theorem 1.6 and Corollary 1.8]{HK}):
\vskip.2cm

 \cite{HK} {\it  Let $A: X \to GL(n,\R)$ be a measurable cocycle over
 an ergodic transformation $T$ of a measure space $(X,\mu)$.
 Then $A$ is measurably cohomologous to a cocycle $B$ with 
 values in an amenable subgroup of $GL(n,\R)$.}

\begin{corollary} \label{dichotomy} 
Let $f$ be a diffeomorphism of a smooth manifold $\M$ preserving
an ergodic measure $\mu$ and let $E$ be a 2-dimensional measurable
invariant distribution. Then the  restriction of the derivative $df$ to $E$ 
has either a measurable invariant one-dimensional distribution or 
a measurable invariant conformal structure.
\end{corollary}

\proof
We apply the Amenable Reduction Theorem to restriction of the 
derivative cocycle to $E$. In the case of $n=2$ any maximal 
amenable subgroup of $GL(2,\R)$ is conjugate either to the
subgroup of (upper) triangular matrices or to the subgroup of
scalar multiples of orthogonal matrices. The former subgroup
preserves a coordinate line, the latter subgroup preserves the 
standard conformal structure on $\R^2$. The measurable coordinate 
change given by the cohomology gives the corresponding 
measurable invariant structure for $df|_E$.
$\QED$
\vskip.2cm

We will now show that $f$ is uniformly quasiconformal on $E^u$ and $E^s$.  
The corollary above applied to $E=E^u$
yields either a measurable invariant one-dimensional distribution 
in $E^u$ or a measurable invariant conformal structure on $E^u$.
In the latter case Proposition \ref{structure} implies that the conformal 
structure is continuous. If follows that  $f$ is  uniformly quasiconformal 
on $E^u$.

If there is a measurable invariant one-dimensional distribution 
in $E^u$ then it is continuous by Proposition \ref{distribution}. For any 
periodic point $p$ this gives an invariant line for the derivative 
$df^n |_{E^u (p)}$ of the return map. This implies that the eigenvalues 
of $df^n |_{E^u (p)}$ are real. Hence, by the assumption of the theorem,
they are either $\lambda ,\lambda$ or $\lambda ,-\lambda$. In the former 
case we have $df^n |_{E^u (p)} = \lambda \cdot \Id$ which preserves any 
conformal structure on $E^u(p)$. This implies uniform quasiconformality 
of $f$ on $E^u$  (see \cite[Proposition 1.1]{KaS}, \cite{L3}).
The case $\lambda ,-\lambda$ is orientation reversing and can be excluded 
by passing to a finite cover of $\M$ and replacing $f$ by its power. Note that
uniform quasiconformality of a power of $f$ implies that of $f$.

We conclude that $f$ is uniformly quasiconformal on $E^u$.
Similarly, $f$ is uniformly quasiconformal on $E^s$. This 
implies that a finite cover of $f$ is $C^\infty$ conjugate to an Anosov 
automorphism of a torus. Indeed, if $\M$ is known to be an infranilmanifold
then \cite[Theorem 1.1]{KaS} stated in the introduction applies. 
Since $f$ is volume preserving the result holds even if $\M$ is 
arbitrary \cite[Corollary 1]{F}. 
$\QED$

\subsection{Proof of Corollary \ref{perturbation}}

The coincidence of the periodic data implies that Theorem \ref{dim4} 
applies to $g$ as well. Hence, after passing to a finite cover of $\M$,  
$f$ and $g$ are $C^\infty$ conjugate to automorphisms of $\T^4$. 
Thus it suffices to show that any two Anosov 
automorphisms $A$ and $B$ of $\T^k$ which are topologically conjugate 
are also $C^\infty$ conjugate. Let $h$ be a topological conjugacy, i.e. a
homeomorphism of $\T^k$ such that $A\circ h= h\circ B$. 
Let $H$ be the induced action of $h$ on the fundamental group $\Z^k$ 
of  $\T^k$. Then $H$ is an integral matrix with determinant $\pm 1$,
and hence it induces an automorphism of $\T^k$.
From the induced actions of $A$, $B$, and $h$ on the fundamental group 
$\Z^k$ we see that $A\circ H= H\circ B$. Thus, $H$ gives a smooth conjugacy 
between $A$ and $B$.
In fact, $H=h$ since the conjugacy to an Anosov automorphism is known 
to be unique in a given homotopy class \cite{KH}.  
$\QED$

\subsection{Proof of Corollary \ref{dim3}}
Let $h$ be the topological conjugacy between $f$ and $g$ given by the 
structural stability. Suppose that $f$ and $g$ have two-dimensional unstable
and one-dimensional stable distributions. We apply the approaches used by 
R. de la Llave in \cite{L0,L1} for one-dimensional distributions and in \cite{L4} 
for higher-dimensional conformal case. There the smoothness of $h$ is established 
separately along the stable and the unstable foliations. This implies that $h$ is 
smooth on the manifold.

The smoothness of $h$ along the two-dimensional unstable foliation follows as in the
proof of Theorem 1.1 in \cite{L4} once we establish that $f$ and $g$ preserve $C^1$
conformal structures on their unstable distributions. As in the proof of Theorem 
\ref{dim4}, we obtain that $f$ is uniformly quasiconformal on $E^u$, and hence
it preserves a bounded measurable conformal structure on $E^u$. 
By (4) of Theorem \ref{main}, this conformal structure is $C^1$. The same argument 
gives an invariant $C^1$ conformal structure for $g$. We conclude that $h$ is smooth
along the unstable foliations. The smoothness along the one-dimensional stable 
foliation can be established as in \cite{L0,L1}.
$\QED$


\end{document}